\documentclass{amsart}

\usepackage{amssymb,amsmath,amsthm,color}

\usepackage[bookmarksopen=true,bookmarksnumbered=true, bookmarksopenlevel=2, colorlinks=true,linkcolor=mblue,
citecolor=mblue,urlcolor=mblue, pdfstartview=FitH]{hyperref}

\definecolor{mblue}{rgb}{0,0,.8}

\usepackage[all]{xy}

\newcommand{\N}{\mathbb N}
\newcommand{\Z}{\mathbb Z}
\newcommand{\Q}{\mathbb Q}
\newcommand{\R}{\mathbb R}
\newcommand{\C}{\mathbb C}
\newcommand{\T}{\mathbb T}
\newcommand{\D}{{\mathcal D}}

\newtheorem{thm}{Theorem}
\newtheorem{lem}{Lemma}
\newtheorem{dfn}{Definition}
\newtheorem{prop}{Proposition}

\begin{document}
\title[Hecke operators in spaces of Maass wave forms]{New models for the action of
Hecke operators in spaces of Maass wave forms.}
\author{Ian Kiming}
\address{Dept. of Mathematics, University of Copenhagen, Universitetsparken 5, 2100 Copenhagen \O , Denmark}

\begin{abstract} Utilizing the theory of the Poisson transform, we develop some new concrete models for the Hecke theory in a
space $M_{\lambda}(N)$ of Maass forms with eigenvalue $1/4-\lambda^2$ on a congruence subgroup $\Gamma_1(N)$. We
introduce the field $F_{\lambda} = {\mathbb Q} (\lambda ,\sqrt{n} , n^{\lambda /2} \mid ~n\in \N )$ so that
$F_{\lambda}$ consists entirely of algebraic numbers if $\lambda = 0$.

The main result of the paper is the following. For a packet $\Phi = (\nu_p \mid p\nmid N)$ of Hecke eigenvalues
occurring in $M_{\lambda}(N)$ we then have that either every $\nu_p$ is algebraic over $F_{\lambda}$, or else $\Phi$
will - for some $m\in {\N}$ - occur in the first cohomology of a certain space $W_{\lambda,m}$ which is a space of
continuous functions on the unit circle with an action of $\mathrm{SL}_2({\R})$ well-known from the theory of
(non-unitary) principal representations of $\mathrm{SL}_2({\R})$.
\end{abstract}

\maketitle

\section{Introduction} \label{intro} We shall fix the following
notation throughout the paper: $N$ is a natural number, $p$ always denotes a prime number, $\Gamma$ is the congruence
subgroup $\Gamma_1(N) \le \mathrm{SL}_2(\Z )$, $G:=\mathrm{SL}_2(\R )$, ${\mathbf T} :=\R /\pi \Z$, $\lambda$ is a
complex number, and $M_{\lambda}(N)$ denotes the complex vector space of Maass forms $g$ on $\Gamma$ with eigenvalue
$\frac{1}{4} -\lambda ^2$ for the Laplacian (\cite{M}), i.e. real-analytic functions $g(x,y)$ on the upper half plane
$\mathfrak{H} := \{ x+iy \mid ~y>0\}$ with the following 3 properties:
$$
-y^2 \left( \frac{\partial ^2}{\partial x^2} + \frac{\partial ^2}{\partial y^2} \right) g = \left( \frac{1}{4} -\lambda
^2\right) g~, \leqno{(i)}
$$
$$
g(\gamma .z) = g(z),~~z=x+iy,\mbox{ for all } \gamma \in \Gamma ,\leqno{(ii)}
$$
where $\gamma .$ denotes the usual action of $\gamma \in \Gamma$ on $\mathfrak{H}$,

\noindent ({\it iii}) There is a positive number $d$ such that
$$
g(x+iy) = O(y^d) \mbox{  for  } y\rightarrow +\infty
$$ and
$$
g(x+iy) = O(y^{-d}) \mbox{  for  } y\rightarrow 0+,
$$
uniformly in $x$.

As usual, an element $g\in M_{\lambda}(N)$ is said to be a cusp form, if it is bounded as a function on $\mathfrak{H}$.

Let $\Delta_1(N)$ be the set of matrices $\left( \begin{array}{cc} a & b \\ c & d \end{array} \right)$ with $a,b,c,d\in
\Z$, $a\equiv 1~(N)$, $c\equiv 0~(N)$, and positive determinant. Then for $\alpha \in \Delta_1(N)$ we have a Hecke
operator $T_{\alpha}$ acting on $M_{\lambda}(N)$: Suppose that $\Gamma \alpha \Gamma = \bigcup _i \Gamma \alpha_i$ as a
disjoint union, and that $g\in M_{\lambda}(N)$. Then $g\mid T_{\alpha}$ is the form given by:
$$
(g\mid T_{\alpha})(z) := (\det \alpha)^{-1/2} \sum_i g(\alpha_i .z) ~.
$$
For primes $p\nmid N$ we shall denote as usual by $T_p$ the Hecke operator belonging to $\alpha=\left(
\begin{array}{cc} 1 & 0 \\ 0 & p \end{array} \right)$. Also, $\T =\T _N$ will denote the Hecke
algebra generated over $\C$ by the operators $T_p$, $p\nmid N$, acting on $M_{\lambda}(N)$. Thus, $\T$ is a commutative
$\C$-algebra.

If $W$ is a complex vector space which is also a $\T$-module, and if $\Phi = (\nu_p \mid p\nmid N)$ is a system of
complex numbers, we say that $\Phi$ occurs in $W$, if $\T$ has a (non-zero) eigenvector in $W$ such that $\nu_p$ equals
the eigenvalue corresponding to $T_p$, for all $p\nmid N$.

From an arithmetical point of view, the interest in the Hecke theory of the spaces $M_{\lambda}(N)$ centers around the
case $\lambda =0$, because a part of `Langlands' philosophy' predicts that a system $\Phi = (\nu_p \mid p\nmid N)$
occurs as a system of Hecke eigenvalues in the subspace of cusp forms in $M_0(N)$, precisely if there exists an
irreducible Galois representation $\rho :~ \mathrm{Gal} (\bar{\Q} /\Q ) \rightarrow \mathrm{GL} _2(\C )$ with finite
image, Artin conductor dividing $N$, $(\det \rho )(\mbox{complex conjugation} )=1$, and such that $\nu_p = \mathrm{Tr}
\rho (\mbox{Frobenius at } p)$, $\forall p\nmid N$. A consequence of this prediction is the conjecture that for every
such system $\Phi$, every $\nu_p$ is an algebraic integer. In \cite{BCR}, \cite{BCR1}, an interesting attempt was made
to prove this conjecture; however, as discovered by Henniart \cite{Hen}, the approach was unfortunately irreparably
flawed. Thus, the conjecture is - in its full generality - still completely open.
\smallskip

It seems clear that any progress in connection with this problem must involve the development of new and more
structured models for the Hecke theory in the spaces $M_{\lambda}(N)$. The purpose of the present paper is to suggest
and initiate the study of some new such models. Here is a rough description of the ideas involved. First we utilize the
Poisson transform to transform the space $M_{\lambda}(N)$ to a space of distributions on the unit circle ${\mathbf T}$,
- equipped with a certain on $\lambda$ depending action of $\mathrm{GL}^+_2({\R})$. The basic idea is then - roughly
speaking - to study `point values' - and the action of Hecke operators on them - of these distributions for points in a
certain countable dense subset $\Xi \subseteq {\mathbf T}$. The actual realization of this idea is a bit technical: By
representing the space of distributions in question as derivatives of continuous functions, one ultimately ends up -
for some $k\in {\N}_0$ - with a very short exact sequence of (cohomological) Hecke modules
$$
H^0(\Gamma , U_{\lambda ,k}) \longrightarrow M_{\lambda}(N) \longrightarrow H^1(\Gamma ,V_{\lambda ,k})
~;\leqno{(\ast)}
$$
here, $U_{\lambda ,k}$ is the space of $(k+1)$-tuples of continuous functions on ${\mathbf T}$ with absolutely
converging Fourier series and equipped with a certain $\mathrm{GL}^+_2({\R})$-action, and $V_{\lambda ,k}$ is a
subspace of $U_{\lambda ,k}$ that is shown to have a natural filtration whose successive quotients are isomorphic to
spaces $W_{\lambda, m}$, $1\le m\le k$, introduced in definition \ref{d3} below: These are spaces consisting of
continuous functions on ${\mathbf T}$ with a certain growth condition on their Fourier coefficients and equipped with
an - on $(\lambda ,m)$ depending - action of $\mathrm{GL}^+_2({\R})$ that is easily recognizable from the theory of
(non-unitary) principal series representations (of $\mathrm{SL}_2({\R})$). Apart from some remarks at the end of the
paper, we shall make no further references to representation theory of $\mathrm{SL}_2({\R})$ as we shall introduce all
concepts in a completely explicit and self-contained manner.
\smallskip

Using $(\ast)$ combined with a study of the evaluation of elements of $H^0(\Gamma , U_{\lambda ,k})$ at points of the
above $\Xi\subseteq {\mathbf T}$, and by using the Ash-Stevens `lifting lemma' for packages of Hecke eigenvalues
(\cite{AS}), we are then able to prove the following theorem.

\begin{thm} \label{t1} Let $\Phi = (\nu_p \mid p\nmid N)$ be a system of
Hecke eigenvalues occurring in $M_{\lambda}(N)$. Then either

$(1)$ Every $\nu_p$ is algebraic over the field $\Q (\lambda ,\sqrt{n} , n^{\lambda /2} \mid ~n\in \N )$,

\noindent or

$(2)$ $\Phi$ occurs in $H^1(\Gamma ,W_{\lambda ,m})$ for some $m\in \N$.
\end{thm}

\section{} \label{1} In this section we prove Theorem \ref{t1}. We
proceed first by recalling some facts about the action of Hecke operators on cohomology of Hecke-modules, then with
various preliminary constructions before introducing the spaces $W_{\lambda ,m}$ and proving Theorem \ref{t1}.

\subsection{} \label{2.1} In the terminology of \cite{AS}, section 1.1, we will be concerned
with the Hecke algebra $\mathcal H$ of the Hecke pair $(\Gamma ,\Delta_1(N))$, so that our ${\T}$ introduced above is a
subalgebra of ${\mathcal H}\otimes {\C}$. So, if $M$ is a right ${\C}\Delta_1(N)$-module, we have a natural right
action of ${\T}$ on the cohomology groups $H^i(\Gamma , M)$, $i\ge 0$, cf. {\it loc. cit.}. In particular, we have a
right action of ${\T}$ on the $\Gamma$-fixed points of $M$ that we consistently denote $M^{\Gamma}$, and also on
$H^1(\Gamma ,M)$: Explicitly, if $\alpha \in \Delta_1(N)$ with $\Gamma \alpha \Gamma = \cup_i \Gamma \alpha_i$
(disjoint), the action of the Hecke operator $T_{\alpha}$ on a homogeneous 1-cocycle $c\colon \Gamma \times \Gamma
\rightarrow M$ is given by
$$
(c\mid T_{\alpha})(\gamma_0 ,\gamma_1) := \sum_i c(t_i(\gamma_0 ),t_i(\gamma_1 )) \mid \alpha_i ~,
$$
where $t_i\colon \Gamma \rightarrow \Gamma$ is the map determined by the requirements $\Gamma \alpha_i \gamma = \Gamma
\alpha_j$, for some $j$ depending on $i$ and $\gamma$, and $\alpha_i \gamma = t_i(\gamma) \alpha_j$.

Let us also recall that the right action of $\T$ commutes with the long exact cohomology sequence associated with a
short exact sequence of ${\C}\Delta_1(N)$-modules.

\subsection{} \label{2.2} In this subsection, we shall consider the
Poisson-transformation associated with the upper half plane $\mathfrak{H}$ interpreted as the symmetric space
$\mathrm{SL}_2(\R ) / \mathrm{SO}_2(\R )$, use this to pull back the action of $\mathrm{GL}_2^+(\R )$ on
$\Delta$-eigenfunctions on $\mathfrak{H}$ to the space of distributions on ${\mathbf T} =\R /\pi \Z$, and establish an
isomorphism as Hecke-modules between $M_{\lambda}(N)$ and a certain space of distributions on ${\mathbf T}$. Let us
first introduce the appropriate actions:
\smallskip

For $\gamma = \left( \begin{array}{cc} a & b \\ c & d \end{array} \right) \in \mathrm{GL}_2^+(\R )$ and $\theta \in
{\mathbf T}$, we put:
$$
j(\gamma , \theta ) := \left( (a\cos \theta +b\sin \theta )^2 + (c\cos \theta + d\sin \theta )^2 \right) ^{1/2} ~,
$$
so that for fixed $\gamma$, $j(\gamma ,\cdot )$ is a $C^{\infty}$ function on ${\mathbf T}$. Secondly, we define
$\gamma .\theta \in {\mathbf T}$ by the requirement
$$
(\cos \gamma .\theta , \sin \gamma .\theta ) = \pm \left( \frac{a\cos \theta +b\sin \theta}{j(\gamma , \theta )} ,
\frac{c\cos \theta + d\sin \theta}{j(\gamma , \theta )} \right) ~,
$$
where the sign is chosen such that $\gamma .\theta \in [0,\pi [$. One immediately verifies that $(\gamma ,\theta )
\mapsto \gamma .\theta$ actually defines an action of $\mathrm{GL}_2^+(\R )$ on ${\mathbf T}$, and that we have
$$
j(\gamma_1 \gamma_2 ,\theta ) = j(\gamma_1 ,\gamma_2 .\theta )j(\gamma_2 ,\theta ) ~.
$$
With this, we see that we have a right action $\mid _{\lambda}$ of $\mathrm{GL}_2^+(\R )$ on $C^{\infty}$-functions on
the torus ${\mathbf T}$ given by:
$$
(\varphi \mid _{\lambda} \gamma)(\theta) := (\det \gamma)^{1+\lambda /2} j(\gamma,\theta)^{-1-\lambda} \cdot \varphi
(\gamma .\theta)~,\quad \mbox{for } \varphi \in C^{\infty}({\mathbf T}) ~.
$$

\begin{dfn} \label{d1} Define ${\D}_{\lambda}$ to be the complex vector space
of distributions on ${\mathbf T}$ with the (on $\lambda$ depending) right action of $\mathrm{GL}_2^+(\R )$ given
by:
$$
(\Lambda \mid _{\lambda} \gamma) \varphi := \Lambda (\varphi \mid _{\lambda} \gamma^{-1} )
$$
for $\Lambda \in {\D}_{\lambda}$, $\varphi \in C^{\infty}({\mathbf T})$, and $\gamma \in \mathrm{GL}_2^+(\R )$.
\end{dfn}

With this definition, the space ${\D}_{\lambda}^{\Gamma}$ is also endowed with the structure of a $\T$-module.
Explicitly, if $\Lambda \in {\D}_{\lambda}^{\Gamma}$ and $\alpha \in \Delta_1(N)$ with $\Gamma \alpha \Gamma = \bigcup
_i \Gamma \alpha_i$ as a disjoint union, then
$$
\Lambda \mid T_{\alpha} := \sum_i (\Lambda \mid _{\lambda} \alpha_i) ~.
$$

\begin{prop} \label{p1} Suppose that $\mathrm{Re}(\lambda) \ge 0$. Then
$$
{\D}_{\lambda}^{\Gamma} \cong M_{\lambda}(N)
$$
as $\T$-modules.
\end{prop}

\begin{proof} Consider the maximal compact subgroup
$$
K:=\mathrm{SO}_2(\R ) = \left\{ r(\theta) := \left(
\begin{array}{cc} \cos \theta & -\sin \theta \\ \sin \theta & \cos \theta \end{array} \right)
\mid ~\theta \in ]-\pi , \pi] \right\}
$$
in $G$. We identify $\mathfrak{H}$ with the symmetric space $G/K$; explicitly,
$x+iy\in \mathfrak{H}$ is identified with the coset $g_{x,y}K$ where
$$
g_{x,y} := \left( \begin{array}{cc} y^{1/2} & xy^{-1/2} \\ 0 & y^{-1/2} \end{array} \right) ~.
$$
We shall now utilize the Poisson transform and the Helgason isomorphism associated with this situation. We shall use
the particular version given in \cite{fj}, Chap. IV, Theorem 5, and proceed now with introducing the necessary
notation.
\smallskip

Consider the following standard subgroups of $G$:
$$
A:=\left\{ \left( \begin{array}{cc} a & 0 \\ 0 & a^{-1} \end{array} \right) \mid ~a\in {\R}_+ \right\} ~,
$$
whose Lie algebra $\mathfrak{a}$ we identify with $\R$, and
$$
N:=\left\{ \left( \begin{array}{cc} 1 & t \\ 0 & 1
\end{array} \right) \mid ~t\in {\R} \right\} ,\quad M:=\left\{ \left( \begin{array}{cc} \pm 1 & 0 \\ 0 & \pm 1
\end{array} \right) \right\} ,
$$
so that we have the Iwasawa decomposition $G=KAN$, and can identify $K/M$ with ${\mathbf T}$. Explicitly, $r(\theta)M$
is identified with
$$
\theta \pmod \pi \in {\mathbf T} ~,
$$
and the Iwasawa decomposition of an element $g=\left( \begin{array}{cc} a & b \\ c & d
\end{array} \right) \in G$ has the shape
$$
g=r(\theta) \left( \begin{array}{cc} u & 0 \\ 0 & u^{-1} \end{array} \right) \left(
\begin{array}{cc} 1 & v \\ 0 & 1 \end{array} \right) ~,
$$
with $u:=(a^2+c^2)^{1/2}$ and $\theta \in ]-\pi ,\pi ]$ determined by $\cos \theta = a/u$, $\sin \theta = c/u$. The map
$h\colon G\longrightarrow \mathfrak{a}$ is the map uniquely determined by the requirement
$$
g \in K\exp (h(g))N, \quad \mbox{for } g\in G ~.
$$
Explicitly, one finds (with the identification of $\mathfrak{a}$ with $\R$)
$$
h\left( \left( \begin{array}{cc} a & b \\ c & d \end{array} \right) \right) = \log (a^2+c^2)^{1/2} ~.
$$
Finally, define for $g\in G$:
$$
|g| := \mathrm{trace}_{\mathfrak{g}} (\mathrm{ad}(gg^t))~,
$$
where $g^t$ is the transpose of $g$, and $\mathfrak{g}$ the Lie algebra of $G$. One shows that $|k_1gk_2| = |g|$ for
$k_1,k_2\in K$, so that it makes sense to define the notion of {\it at the most exponential growth} for functions $f$
on $G/K$ by the requirement
$$
|f(g)| \le a|g|^b
$$
with some constants $a\in {\R}_+$, $b\in \R$.
\medskip

Now we specialize Chap. IV, 2 of \cite{fj}, in particular Theorem 5 (ii), to our present situation: We identify
elements $\lambda$ of the space of linear forms on the complexified Lie algebra $\mathfrak{a}_{\C}$ with complex
numbers. The Poisson transform $\mathcal{P}_{\lambda}$ associated with $\lambda$ is then a linear map from the space
${\D}(K/M)$ of distributions on $K/M$ to the space $E_{\lambda}$ of $C^{\infty}$ functions on $G/K$ with at the most
exponential growth and are eigenfunctions for the Laplace-Beltrami operator on $G/K$ with an eigenvalue that we specify
below; $\mathcal{P}_{\lambda}$ is given explicitly by
$$
(\mathcal{P}_{\lambda} \Lambda )(gK) := \Lambda \left( kM \mapsto e^{(-\lambda -1)\cdot h(g^{-1}k)} ~,~k\in K \right)
~,\leqno{(\ast)}
$$
for $\Lambda \in {\D}(K/M)$. Furthermore, Theorem 5 (ii) of {\it loc. cit.} implies that $\mathcal{P}_{\lambda}$ is an
{\it isomorphism}, provided that $\mathrm{Re}(\lambda) \ge 0$ which is insured in the present case by assumption. Using
the above identification of $K/M$ with $\mathbf T$ and of $G/K$ with $\mathfrak{H}$ via $x+iy\leftrightarrow g_{x,y}K$,
it is straightforward to verify that we can rewrite $(\ast)$ as
$$
(\mathcal{P}_{\lambda} \Lambda )(x+iy) = \Lambda \left( \theta \mapsto b(x,y,\theta)^{\frac{-1-\lambda}{2}} ~,~\theta
\in {\mathbf T} \right) ,\leqno{(\ast \ast)}
$$
where
$$
b(x,y,\theta) := y^{-1} \left( \cos ^2\theta + x^2\sin ^2\theta - 2x\cos \theta \sin \theta \right) + y\sin ^2\theta
~,
$$
which arises because of
$$
h(g_{x,y}^{-1} r(\theta)) = \frac{1}{2} \log b(x,y,\theta) ~.
$$
From this, one readily checks that elements of $E_{\lambda}$ viewed as $C^{\infty}$ functions on $\mathfrak{H}$ have
$\Delta$-eigenvalue $1/4-\lambda^2$. Furthermore, it is tedious but straightforward to check that we have defined the
action of $\mathrm{GL}_2^+({\R})$ on $\D$ precisely so as to make $\mathcal{P}_{\lambda}$ equivariant w.r.t. to
$\mathrm{GL}_2^+({\R})$-action. Thus, we may conclude that
$$
{\D}^{\Gamma} \cong E_{\lambda}^{\Gamma} ~,
$$
as $\T$-modules. Hence the proof is concluded by showing that
$$
M_{\lambda}(N) = E_{\lambda}^{\Gamma} ~.
$$
Now, the growth condition on elements of $f\in E_{\lambda}$ viewed as functions on $\mathfrak{H}$ requires the
existence of constants $a\in {\R}_+$ and $b\in \R$ (depending on $f$) such that
$$
|f(x+iy)| \le a|g_{x,y}|^b ~,
$$
and as a simple computation shows that
$$
|g_{x,y}| = 1 + 2 x^2y^{-2} + (y + x^2y^{-1} )^2 + y^{-2} \ge \max \{ y^2,y^{-2} \} ~,
$$
it is immediately clear that $M_{\lambda}(N) \le E_{\lambda}^{\Gamma}$. On the other hand, one sees that $|g_{x,y}| \le
(\mathrm{const.}) \cdot (y^2+y^{-2})$ for $x\in [0,1]$, so that if $f\in E_{\lambda}^{\Gamma}$ then certainly for some
positive $d$ we have
$$
f(x+iy) = O(y^{d}) \quad \mbox{as  } y\rightarrow +\infty
$$
and
$$
f(x+iy) = O(y^{-d}) \quad \mbox{as  } y\rightarrow 0+
$$
uniformly in $x\in [0,1]$. But since $f$ is invariant under the substitution $x\mapsto x+1$, this holds uniformly in
$x\in \R$. Thus, $f\in M_{\lambda}(N)$.
\end{proof}

If $f$ is a continuous complex-valued function on ${\mathbf T}$, we shall write
$$
f \sim \sum_{n\in \Z} a_n e^{2in\theta}
$$
(only) to signal that the Fourier coefficients of $f$ are being denoted by $a_n$, $n\in \Z$.
\smallskip

If $\Lambda \in {\D}({\mathbf T})$, the Fourier coefficients of $\Lambda$ are
$$
a_n:= \Lambda (e^{-2in\theta}),
$$
and as is well-known, a given sequence $(a_n)$ of complex numbers is the sequence of Fourier coefficients of some
distribution $\Lambda \in {\D}({\mathbf T})$ if and only if there exists an integer $m\ge 0$ such that
$$
\sum_n (1+n^2)^{-m} |a_n|^2 < \infty ~;
$$
if this is the case, there is then an integer $s \ge 0$ such that the series
$$
\sum_n (1+2in)^{-s} \cdot a_n
$$
converges absolutely. Then
$$
f(\theta) := \sum_n (1+2in)^{-s} a_n e^{2in\theta}
$$
defines a continuous function on $\mathbf T$, and with the notation
$$
\partial := 1+ \frac{d}{d\theta} ~,
$$
we have $\Lambda = \partial ^s f$ in the usual sense, i.e.
$$
\Lambda \varphi = \frac{1}{\pi} \int_{\mathbf T} f(\theta) (\partial ^s\varphi )(\theta) d\theta
$$
for $\varphi \in C^{\infty}({\mathbf T})$. It will be convenient for us to consider distributions representable in
slightly more general form
$$
\Lambda = \sum_{j=0}^s \partial ^j f_j~,\leqno{(\dag)}
$$
where the $f_j$ are continuous functions on $\mathbf T$ {\it with absolutely converging Fourier series}. We now define
some auxiliary objects, and then the spaces occurring in Theorem \ref{t1}.

\begin{dfn} \label{d2} For $s \in {\N}_0$, denote by ${\D}_{\lambda ,s}$ the subspace of ${\D}({\mathbf T})$
consisting of distributions $\Lambda$ representable in the form $(\dag)$ and with the right
$\mathrm{GL}_2^+({\R})$-action given in Definition \ref{d1}. Define also $U_{\lambda ,s}$ to be the complex vector
space consisting of $(s+1)$-tuples $(f_0,\ldots ,f_s)$ of continuous functions on $\mathbf T$ with absolutely
converging Fourier series, and define $V_{\lambda ,s}$ to be the subspace of $U_{\lambda ,s}$ consisting of tuples
$(f_0,\ldots ,f_s)$ with
$$
\sum_j \partial ^j f_j = 0
$$
in the distribution sense. Because of Proposition \ref{p1} and \cite{M}, Satz 5, the space ${\D}_{\lambda}^{\Gamma}$ is
finite-dimensional. Consequently, there exists a non-negative integer $k$ such that
$$
{\D}_{\lambda}^{\Gamma} \le {\D}_{\lambda ,k} ~,
$$
so that obviously,
$$
{\D}_{\lambda ,k}^{\Gamma} = {\D}_{\lambda}^{\Gamma} ~.
$$
\end{dfn}

\begin{dfn} \label{d3} Let $m\in \N$. We define $W_{\lambda ,m}$ as the complex vector space consisting of
continuous functions $f \sim \sum_{n\in \Z} a_n e^{2in\theta}$ on ${\mathbf T}$ with
$$
\sum_{n\in \N} |n| |a_n| < \infty ~,
$$
and with the following on $(\lambda ,m)$ depending right action of $\mathrm{GL}_2^+({\R})$:
$$
(f\mid _{\lambda ,m} \gamma )(\theta) := f(\gamma .\theta ) (\det \gamma )^{-m-\lambda /2} j(\gamma ,\theta
)^{(2m-1)+\lambda} ~.
$$
\end{dfn} Using the fact that $j(\gamma ,\cdot )\in C^{\infty} ({\mathbf T})$, we have
$j(\gamma ,\cdot )\in W_{\lambda ,m}$ for every $m\in \N$, and it is easily verified that $\cdot \mid _{\lambda ,m}
\gamma$ actually maps $W_{\lambda ,m}$ into itself so that the definition makes sense.
\smallskip

The definitions of $U_{\lambda ,s}$ and $V_{\lambda ,s}$ are such that we have a natural exact sequence of complex
vector spaces
$$
0\longrightarrow V_{\lambda ,s} \longrightarrow U_{\lambda ,s} \longrightarrow {\D}_{\lambda ,s} \longrightarrow 0
~,
$$
where the map $U_{\lambda ,s} \longrightarrow {\D}_{\lambda ,s}$ is given by $(f_0,\ldots ,f_s) \mapsto \sum_j
\partial ^j f_j$. In the next subsection we lift the action of $\mathrm{GL}_2^+({\R})$ on ${\D}_{\lambda ,s}$ to the
space $U_{\lambda ,s}$ so that this sequence becomes an exact sequence of $\Delta_1(N)$-modules, and we show that the
$W_{\lambda ,m}$ appear as subquotients of the modules $V_{\lambda ,s}$. In subsection \ref{2.4}, we prove a theorem
concerning eigenvalues of Hecke operators acting on the space $U_{\lambda ,s}^{\Gamma}$ of $\Gamma$-fixed points of
$U_{\lambda ,s}$. With this preparation, we then proceed in subsection \ref{2.5} to finish the proof of Theorem
\ref{t1}.

\subsection{} \label{2.3}

\begin{lem} \label{l1} For each $s\in {\N}_0$ and to each $\gamma \in \mathrm{GL}_2^+({\R})$ there exist on $\lambda$
depending, uniquely determined $C^{\infty}$ functions $u_{\gamma}^{s,j}(\theta)$, $j\in {\Z}$, on $\mathbf T$ with
$u_{\gamma}^{s,j} = 0$ for $j<0$ and for $j>s$, and with the following properties.
\smallskip

\noindent $(1)$  For all $\varphi \in C^{\infty}({\mathbf T})$ we have
$$
\partial ^s(\varphi \mid _{\lambda} \gamma )(\theta) =
(\det \gamma )^{1+\lambda /2} \sum_{j=0}^s u_{\gamma}^{s,j}(\theta) \cdot (\partial ^j\varphi )(\gamma .\theta ) ~.
$$

\noindent $(2)$  For all $s$,
$$
u_{\gamma}^{s,s}(\theta) = j(\gamma ,\theta )^{-1-\lambda -2s} \cdot (\det \gamma)^s .
$$

\noindent $(3)$  For all $s$ and $j=0,\ldots s$, the function $j(\gamma ,\theta)^{\lambda} u_{\gamma}^{s,j}(\theta)$ is a
polynomial in the functions
$$
\frac{d^{\ell}}{d\theta ^{\ell}} j(\gamma ,\theta) \quad , ~\ell =0,\ldots ,s-j,
$$
with coefficients in ${\Z}[\lambda ,(\det \gamma)]$.
\smallskip

\noindent $(4)$  For $\gamma_1 ,\gamma_2 \in \mathrm{GL}_2^+({\R})$,
$$
u_{\gamma_1 \gamma_2}^{s,j}(\theta) = \sum_{\ell =j}^s u_{\gamma_2}^{s,\ell}(\theta) \cdot u_{\gamma_1}^{\ell
,j}(\gamma_2 .\theta) ~.
$$
\end{lem}

\begin{proof} This is a completely trivial exercise in differential calculus so we shall be very brief. Define first
$u_{\gamma}^{s,j}(\theta) := 0$ for $j<0$ and for $j>s$. Put $u_{\gamma}^{0,0}(\theta) := j(\gamma ,\theta )^{-1-\lambda}$, and
recursively (w.r.t. $s$)
$$
u_{\gamma}^{s+1,j}(\theta) := \partial u_{\gamma}^{s,j}(\theta) + (\det \gamma) \cdot j(\gamma
,\theta)^{-2} \left( u_{\gamma}^{s,j-1}(\theta) - u_{\gamma}^{s,j}(\theta) \right) ,
$$
for $j=0,\ldots s+1$. Then (2) and (3) are immediately clear, and (1) is easily proved by induction on $s$ using the
relation $\partial (fg) = -fg + f(\partial g) + (\partial f)g$ and that
$$
\frac{d(\gamma .\theta)}{d\theta} = (\det \gamma) \cdot j(\gamma ,\theta)^{-2} ~.
$$
Notice that the case $s=0$ in (1) is merely the definition of the action $\mid _{\lambda} \gamma$ on $C^{\infty}$
functions. Uniqueness for a given $\gamma$ of functions $u_{\gamma}^{s,j}(\theta)$ with property (1), is shown by
assuming functions $w_{\gamma}^{s,j}(\theta)$ given with
$$
0= \sum_{j=0}^s w_{\gamma}^{s,j}(\theta) \cdot (\partial ^j\varphi )(\gamma .\theta )
$$
for all $\varphi \in C^{\infty}({\mathbf T})$. Using this on test functions $\varphi (\theta) = e^{2im\theta}$, $m\in
{\N}$, changing variables $\theta \mapsto \gamma ^{-1} .\theta$, and letting $m\rightarrow \infty$, one obtains
successively $w_{\gamma}^{s,s} = 0, \ldots , w_{\gamma}^{s,0} = 0$.
\smallskip

Finally, uniqueness for each fixed $\gamma$ of functions $u_{\gamma}^{s,j}$ with (1) proves (4): Applying $\partial ^s$ to the a
function $\varphi \mid _{\lambda} \gamma_1 \gamma_2 = (\varphi \mid _{\lambda} \gamma_1) \mid _{\lambda} \gamma_2$, we see that
(1) holds for $\gamma = \gamma_1 \gamma_2$ if $u_{\gamma_1 \gamma_2}^{s,j}$ is replaced by the function
$$
\sum_{\ell =0}^s u_{\gamma_2}^{s,\ell}(\theta) \cdot u_{\gamma_1}^{\ell ,j}(\gamma_2 .\theta) = \sum_{\ell =j}^s
u_{\gamma_2}^{s,\ell}(\theta) \cdot u_{\gamma_1}^{\ell ,j}(\gamma_2 .\theta) ~.
$$
By uniqueness, (4) follows.
\end{proof}

\begin{dfn} \label{d4} Define for $s\in {\N}_0$ and $\gamma \in \mathrm{GL}_2^+({\R})$ a matrix function (depending
on $\lambda$) of size $(s+1) \times (s+1)$ by
$$
A_{\gamma}^{(s)} (\theta ) := \left( (\det \gamma)^{-\lambda /2} j(\gamma ,\theta)^{-2} u_{\gamma^{-1}}^{\mu ,\nu}
(\gamma .\theta) \right) _{0\le \mu ,\nu \le s} ~,
$$
$\theta \in {\mathbf T}$, where the $ u_{\gamma^{-1}}^{\mu ,\nu}$ are the uniquely determined functions from Lemma
\ref{l1}.

Define also for $(f_0,\ldots ,f_s) \in U_{\lambda ,s}$
$$
(f_0,\ldots ,f_s) \mid _{\lambda} \gamma = (\tilde{f} _0,\ldots ,\tilde{f} _s) ~,
$$
where
$$(\tilde{f} _0(\theta ),\ldots ,\tilde{f} _s(\theta) := (f_0(\gamma .\theta) ,\ldots ,f_s(\gamma .\theta) )
\cdot A_{\gamma}^{(s)} (\theta) ~,
$$
and the linear map $h_s \colon U_{\lambda ,s} \longrightarrow {\D}_{\lambda ,s}$ by
$$
h_s(f_0,\ldots ,f_s) := \sum_j \partial ^j f_j ~.
$$
\end{dfn}

\begin{prop} \label{p2}

\noindent $(1)$  The map $\mid_{\lambda} \gamma$ defines a right $\mathrm{GL}_2^+({\R})$-action on $U_{\lambda ,s}$
such that $h_s$ is a homomorphism of $\mathrm{GL}_2^+({\R})$-modules.
\smallskip

\noindent $(2)$  If $s>0$, the natural injection $\iota_{s-1} \colon U_{\lambda ,s-1} \longrightarrow U_{\lambda ,s}$
given by
$$
(f_0,\ldots ,f_{s-1} ) \mapsto (f_0,\ldots ,f_{s-1},0 )
$$
is a $\mathrm{GL}_2^+({\R})$-homomorphism that injects $V_{\lambda ,s-1}$ into $V_{\lambda ,s}$.
\smallskip

\noindent $(3)$  For $s>0$, if $V_{\lambda ,s-1}$ is viewed as a $\mathrm{GL}_2^+({\R})$-submodule of $V_{\lambda ,s}$,
we have
$$
V_{\lambda ,s}/V_{\lambda ,s-1} \cong W_{\lambda ,s}
$$
as $\mathrm{GL}_2^+({\R})$-modules, and hence also as $\T$-modules.
\end{prop}

\begin{proof} (1) Notice first that $A_{\gamma}^{(s)} (\theta)$ has $C^{\infty}$ coefficients so that
$\mid_{\lambda} \gamma$ actually maps $U_{\lambda ,s}$ into itself. That we have thus defined a right
$\mathrm{GL}_2^+({\R})$-action follows from the formula
$$
A_{\gamma_1 \gamma_2}^{(s)} (\theta) = A_{\gamma_1}^{(s)} (\gamma_2 .\theta) A_{\gamma_2}^{(s)} (\theta)
$$
which is easily verified on the basis of Lemma \ref{l1}, (4).
\smallskip

To see that $h_s$ commutes with the action of $\mathrm{GL}_2^+({\R})$, suppose that
$$
F:=(f_0,\ldots ,f_s) \in U_{\lambda ,s}~,
$$
put $\Lambda := h_s(F)$, and suppose that $\varphi \in C^{\infty}({\mathbf T})$. Then by definition of $\mid_{\lambda}
\gamma$ and $A_{\gamma}^{(s)} (\theta)$, the function $h_s(F\mid_{\lambda} \gamma )\varphi$ equals
$$
\frac{1}{\pi} \int_{\mathbf T} \sum_{\ell ,j} \left( f_j(\gamma .\theta) \cdot (\det \gamma)^{-\lambda /2} j(\gamma
,\theta)^{-2} u_{\gamma^{-1}}^{j,\ell} (\gamma .\theta) \right) \cdot (\partial ^{\ell} \varphi )(\theta) d\theta
,
$$
where the summation is over $\ell ,j =0,\ldots ,s$. Making the change of variables $\theta \mapsto \gamma^{-1}
.\theta$, using $\frac{d(\gamma .\theta)}{d\theta} = (\det \gamma) \cdot j(\gamma ,\theta)^{-2}$, and remembering the
properties of the $u_{\gamma^{-1}}^{j,\ell}$, we find that this integral equals
$$
\Lambda (\varphi \mid_{\lambda} \gamma^{-1} ) = (\Lambda \mid_{\lambda} \gamma) \varphi ~.
$$

\noindent (2)  That $\iota_{s-1}$ is a $\mathrm{GL}_2^+({\R})$-homomorphism is a consequence of the definition of the
action of $\mathrm{GL}_2^+({\R})$ combined with the observation that the matrix $A_{\gamma}^{(s)} (\theta)$ is an upper
triangular $(s+1) \times (s+1)$ matrix whose upper left $s\times s$ minor coincides with $A_{\gamma}^{(s-1)} (\theta)$.
That $\iota_{s-1}$ maps $V_{\lambda ,s-1}$ into $V_{\lambda ,s}$ is trivial.
\smallskip

\noindent (3)  Suppose that $s>0$ and that $(f_0,\ldots ,f_s) \in V_{\lambda ,s}$. Denote the Fourier coefficients of
$f_j$ by $a_n(f_j)$, $n\in {\Z}$. By definition of $V_{\lambda ,s}$, the distribution $\Lambda = \sum_j \partial^j f_j
\in {\D}({\mathbf T})$ is $0$, and hence every Fourier coefficient $a_n(\Lambda)$ of $\Lambda$ vanishes. But we
have
$$
a_n(\Lambda) = \sum_{j=0}^s (1+2in)^j a_n(f_j) ~,
$$
hence
$$
(1+2in) a_n(f_s) = -\sum_{j=0}^{s-1} (1+2in)^{j-s+1} a_n(f_j)
$$
for each $n\in {\Z}$. As each of the series $\sum_n a_n(f_j)$ is absolutely convergent, the same holds for the series
$\sum_n (1+2in) a_n(f_s)$. It follows that $f_s$ lies in the space $W_{\lambda ,s}$. Consequently, we can define a
linear map $\psi_s \colon V_{\lambda ,s} \longrightarrow W_{\lambda ,s}$ by
$$
\psi_s (f_0,\ldots ,f_s) := f_s ~.
$$
Then clearly the kernel of $\psi_s$ is $V_{\lambda ,s-1}$ viewed as a subspace of $V_{\lambda ,s}$ via $\iota_{s-1}$.
We claim that $\psi_s$ is surjective. Suppose that $f\in W_{\lambda ,s}$ with Fourier coefficients $a_n(f)$. Then if we
define the numbers
$$
a_n := -(1+2in)a_n(f)
$$
for $n\in {\Z}$, the series $\sum_n a_n$ converges absolutely. There is thus a continuous function $g$ on ${\mathbf T}$
with Fourier coefficients $a_n(g) = a_n$, and $F:=(0,\ldots ,0,g,f)$ is an element of $U_{\lambda ,s}$. The Fourier
coefficients of the distribution
$$
\partial^{s-1} g +\partial^s f
$$
are all 0, so $F$ is in fact an element of $V_{\lambda ,s}$. But $\psi_s (F) = f$.
\smallskip

It remains to be seen that $\psi_s$ commutes with the action of the group $\mathrm{GL}_2^+({\R})$. Let ${\mathbf
f}:=(f_0,\ldots ,f_s) \in V_{\lambda ,s}$ and let $\gamma \in \mathrm{GL}_2^+({\R})$. By the definition of ${\mathbf
f}\mid_{\lambda} \gamma$ we find
$$
{\mathbf f}\mid_{\lambda} \gamma = (\ldots ,f_s(\gamma .\theta) (\det \gamma )^{-\lambda /2} j(\gamma ,\theta)^{-2}
u_{\gamma^{-1}}^{s,s} (\gamma .\theta) )~.
$$
Together with Lemma \ref{l1} (2) and $j(\gamma^{-1} ,\gamma .\theta) = j(\gamma ,\theta )^{-1}$, this
gives
$$
(\psi_s({\mathbf f}\mid_{\lambda} \gamma ))(\theta) = f_s(\gamma .\theta) (\det \gamma)^{-s-\lambda /2} j(\gamma
,\theta)^{2s-1+\lambda} ~,
$$
which is precisely the definition of $f_s\mid_{\lambda ,s} \gamma$ for $f_s\in W_{\lambda ,s}$.
\end{proof}

\subsection{} \label{2.4} The purpose of this subsection is to prove a statement concerning possible eigenvalues
of Hecke operators acting on the spaces $U_{\lambda ,s}^{\Gamma}$. We proceed with some preparations.
\smallskip

Let us denote by $F_{\lambda}$ the field occurring in Theorem \ref{t1}, (1), i.e.,
$$
F_{\lambda} := {\Q} (\lambda ,\sqrt{n} , n^{\lambda /2} \mid ~n\in \N ) ~.
$$
Denote also by $\Xi$ the subset of ${\mathbf T} = {\R}/\pi {\Z}$ consisting of those $\theta \in {\mathbf T}$ for which
$\cot \theta \in {\mathbb P}^1({\Q})$, i.e.,
$$
\Xi := \{ 0\} \cup \left\{ 0 \not =\theta \in {\mathbf T} \mid \cot \theta \in {\Q} \right\} ~.
$$

\begin{lem} \label{l2} $(1)$  $\Xi$ is dense in ${\mathbf T}$ and stable under the action  on ${\mathbf T}$ of the
group $\mathrm{GL}_2^+({\Q})$.
\smallskip

\noindent $(2)$  Suppose that $\theta_0 \in \Xi$, $\ell \in {\N}_0$, and $\gamma \in \mathrm{GL}_2^+({\Q})$. Then
$$
\frac{d^{\ell}}{d\theta ^{\ell}} j(\gamma ,\theta) _{\mid \theta =\theta_0} \in F_0 ~.
$$
\smallskip

\noindent $(3)$  Suppose that $\theta_0 \in \Xi$, $s \in {\N}_0$, and $\gamma \in \mathrm{GL}_2^+({\Q})$. Then the
matrix
$$
A_{\gamma}^{(s)}(\theta_0)
$$
(definition \ref{d4}) is a $(s+1)\times (s+1)$ matrix with coefficients in $F_{\lambda}$.
\end{lem}

\begin{proof} (1) The first statement is clear, and the second follows from the formula
$$
\cot \gamma .\theta = \frac{a\cot \theta +b}{c\cot \theta +d}
$$
for $\gamma \in \mathrm{GL}_2^+({\R})$.
\smallskip

\noindent (2) Notice first that $\cos \theta ,\sin \theta \in F_0$ if $\theta \in \Xi$. Now, if $\ell =0$, the
statement is clear if $\sin \theta_0 = 0$, and also if $\sin \theta_0 \not = 0$ since then
$$
j(\gamma ,\theta_0) = \pm \sin \theta_0 \left( (a\cot \theta_0 +b)^2 + (c\cot \theta_0 +d)^2 \right) ^{1/2} \in F_0
~.
$$
For $\ell \ge 1$, the statement follows from this by induction on $\ell$ when one shows by induction on $\ell \ge
1$ that
$$
\frac{d^{\ell}}{d\theta ^{\ell}} j(\gamma ,\theta) = j(\gamma ,\theta )^{-\ell} \cdot p_{\gamma ,\ell}(\theta)
~,
$$
where $p_{\gamma ,\ell}(\theta)$ is a polynomial with rational coefficients in $\cos \theta$, $\sin \theta$ and the
$\frac{d^{\mu}}{d\theta ^{\mu}} j(\gamma ,\theta)$, $\mu =0,\ldots ,\ell -1$.
\smallskip

\noindent (3) This follows from (2), the definition of $A_{\lambda}^{(s)}(\theta)$, and from Lemma \ref{l1}, (3).
\end{proof}

Now let $s\in {\N}_0$ and $\alpha \in \Delta_1(N)$. We shall consider the action of the Hecke operator $T_{\alpha}$ on
a space $U_{\lambda ,s}^{\Gamma} = H^0(\Gamma ,U_{\lambda ,s})$. The action of $T_{\alpha}$ on this space is as a
cohomological Hecke operator which we can describe explicitly as follows. Let
$$
\Gamma \alpha \Gamma = \cup_{\mu=1}^r \Gamma \alpha_{\mu} ~,\quad \alpha_{\mu} \in \Delta_1(N) ~,
$$
as a disjoint union. Then if ${\mathbf f} = (f_0,\ldots ,f_s) \in U_{\lambda ,s}^{\Gamma}$, we have
\begin{eqnarray*} ({\mathbf f} \mid_{\lambda} T_{\alpha} ) (\theta)
:= & \sum_{\mu=1}^r ( (f_0,\ldots ,f_s) \mid_{\lambda} \alpha_{\mu} )(\theta) \\
= & \sum_{\mu=1}^r (f_0(\alpha_{\mu} .\theta ),\ldots ,f_s(\alpha_{\mu} .\theta ) )A_{\alpha_{\mu}}^{(s)}(\theta)
.\end{eqnarray*} Let us also as usual denote by $\Gamma_{\infty} :=\langle \pm \left(
\begin{array}{cc} 1 & 1 \\ 0 & 1 \end{array} \right) \rangle$ the
stabilizer in $\Gamma$ of the cusp $\infty \in {\mathbb P}^1({\Q})$. Notice that if $\Gamma$ is viewed as acting on
${\mathbf T}$ then $\Gamma_{\infty}$ is precisely the stabilizer of $0\in {\mathbf T}$. Choose a decomposition of
$\mathrm{SL}_2({\Z})$ in double cosets w.r.t. $(\Gamma ,\Gamma_{\infty})$:
$$
\mathrm{SL}_2({\Z}) = \cup_{\nu=1}^{m} \Gamma \gamma_{\nu} \Gamma_{\infty} ~,
$$
so that the set $\{ \gamma_{\nu} \}$ is in 1-1 correspondence with the cusps w.r.t. $\Gamma$. We define a linear map
$\phi_s$ of $U_{\lambda ,s}^{\Gamma}$ into $\mathrm{Mat}_{m,s+1}({\C}) \cong {\C}^{m(s+1)}$ by
$$
\phi_s (f_0,\ldots ,f_s) := \left( \begin{array}{ccc} f_0(\gamma_1 .0) & \ldots & f_s(\gamma_1 .0) \\
\vdots & & \vdots \\ f_0(\gamma_m .0) & \ldots & f_s(\gamma_m .0) \end{array} \right) ~.
$$

\begin{thm} \label{t2} The space $U_{\lambda ,s}^{\Gamma}$ is finite-dimensional and for any $\alpha \in \Delta_1 (N)$,
the eigenvalues of the Hecke operator $T_{\alpha}$ acting on $U_{\lambda ,s}^{\Gamma}$ are algebraic over the field
$F_{\lambda}$.
\end{thm}

\begin{proof} Retain the above notation. We first show that the linear map $\phi_s$ is injective which will prove the
first part of the theorem. So, suppose that $(f_0,\ldots ,f_s) \in U_{\lambda ,s}^{\Gamma}$ with
$$
\phi_s (f_0,\ldots ,f_s) = 0 ~.
$$
We must show $f_0=\ldots =f_s=0$. View $\mathrm{SL}_2({\Z})$ as acting on ${\mathbf T}$. One finds that
$\mathrm{SL}_2({\Z}) .0 = \Xi$ which is dense in ${\mathbf T}$ (Lemma \ref{l2} (1)). As the $f_j$ are continuous it
thus suffices to show that $f_0(g.0) = \ldots =f_s(g.0) = 0$ for all $g\in \mathrm{SL}_2({\Z})$. Let then $g\in
\mathrm{SL}_2({\Z})$ and write
$$
g=\gamma \cdot \gamma_{\nu} \cdot \gamma_{\infty} ~,
$$
with $1\le \nu \le m$, $\gamma \in \Gamma$, and $\gamma_{\infty} \in \Gamma_{\infty}$. Then,
\begin{eqnarray*}
(f_0(g.0),\ldots ,f_s(g.0)) = & (f_0(\gamma .(\gamma_{\nu} .0)),\ldots ,f_s(\gamma .(\gamma_{\nu} .0)) )
\\ = & ((f_0,\ldots ,f_s) \mid_{\lambda} \gamma )(\gamma_{\nu} .0) \cdot A_{\gamma}^{(s)}(\gamma_{\nu} .0)^{-1} \\ = &
(f_0(\gamma_{\nu} .0),\ldots ,f_s(\gamma_{\nu} .0) )A_{\gamma}^{(s)}(\gamma_{\nu} .0)^{-1}
\\ = & (0,\ldots ,0) .
\end{eqnarray*}
\smallskip

Secondly we show the existence of an endomorphism $t_{\alpha}$ of ${\C}^{m(s+1)}$ with the following 2 properties:
\smallskip

\noindent (i)  The diagram
$$
\xymatrix{ U_{\lambda ,s}^{\Gamma} \ar[d]_{T_{\alpha}} \ar[r]^{\phi_s} & {\C}^{m(s+1)} \ar[d]^{t_{\alpha}} \\
U_{\lambda ,s}^{\Gamma} \ar[r]^{\phi_s}  & {\C}^{m(s+1)} }
$$
is commutative.
\smallskip

\noindent (ii)  $t_{\alpha}$ is defined over the field $F_{\lambda}$ (i.e. given by a matrix with coefficients in
$F_{\lambda}$).
\smallskip

Together with the injectivity of $\phi_s$, this will then prove the rest of the theorem.
\smallskip

Let $\mu \in \{ 1,\ldots ,r\}$ and $\nu \in \{ 1,\ldots ,m\}$. Now, $\alpha_{\mu} \gamma_{\nu} .0 \in \Xi \subseteq
{\mathbf T}$, and as $\Xi = \mathrm{SL}_2({\Z}).0$ there is $g_{\mu ,\nu} \in \mathrm{SL}_2({\Z})$ such that $g_{\mu
,\nu}.0 = \alpha_{\mu} \gamma_{\nu} .0$. We can write
$$
g_{\mu,\nu} = \beta_{\mu,\nu} \gamma_{\xi_{\mu} (\nu)} \gamma_{\infty} ~,
$$
where $\beta_{\mu,\nu} \in \Gamma$, $\xi_{\mu}$ is some map of $\{ 1,\ldots ,m\}$ into itself, and $\gamma_{\infty} \in
\Gamma_{\infty}$. Then
$$
\alpha_{\mu} \gamma_{\nu} .0 = \beta_{\mu,\nu} \gamma_{\xi_{\mu} (\nu)} .0 ~.
$$
Define then the endomorphism $t_{\alpha}$ of $\mathrm{Mat}_{m,s+1}({\C})$ by
$$
t_{\alpha} \left( \begin{array}{ccc} x_{1,0} & \ldots & x_{1,s} \\
\vdots & & \vdots \\ x_{m,0} & \ldots & x_{m,s}
\end{array} \right) :=
$$
$$
\sum_{\mu =1}^r \left( \begin{array}{c} (x_{\xi_{\mu} (1),0}
,\ldots ,x_{\xi_{\mu} (1),s}) A_{\beta_{\mu,1}}^{(s)}(\gamma_{\xi_{\mu} (1)} .0)^{-1} A_{\alpha_{\mu}}^{(s)}(\gamma_1
.0) \\ \vdots \\ (x_{\xi_{\mu} (m),0} ,\ldots ,x_{\xi_{\mu} (m),s}) A_{\beta_{\mu,m}}^{(s)}(\gamma_{\xi_{\mu} (m)}
.0)^{-1} A_{\alpha_{\mu}}^{(s)}(\gamma_m .0) \end{array} \right) ~.
$$
Then claim (ii) above is clear because of Lemma \ref{l2} (3). We proceed to show (i). So, let ${\mathbf f} =
(f_0,\ldots ,f_s) \in U_{\lambda ,s}^{\Gamma}$. Then,
\begin{eqnarray*}
\lefteqn{\phi_s ({\mathbf f} \mid_{\lambda} T_{\alpha} ) = \sum_{\mu =1}^r \phi_s ({\mathbf f} \mid_{\lambda}
\alpha_{\mu} )} \\ = & \sum_{\mu =1}^r \left(
\begin{array}{c} \vdots \\ ( f_0(\alpha_{\mu} \gamma_{\nu}
.0),\ldots ,f_s(\alpha_{\mu} \gamma_{\nu} .0) )
A_{\alpha_{\mu}}^{(s)} (\gamma_{\nu} .0) \\ \vdots \end{array}
\right) \\ = & \sum_{\mu =1}^r \left(
\begin{array}{c} \vdots \\ ( f_0(\beta_{\mu ,\nu}
\gamma_{\xi_{\mu} (\nu)} .0),\ldots ,f_s(\beta_{\mu ,\nu}
\gamma_{\xi_{\mu} (\nu)} .0) ) A_{\alpha_{\mu}}^{(s)}
(\gamma_{\nu} .0) \\ \vdots \end{array} \right) ~.
\end{eqnarray*}

As ${\mathbf f} \in U_{\lambda ,s}^{\Gamma}$ we have
$$
(f_0(\beta_{\mu ,\nu} .\theta ),\ldots ,f_s(\beta_{\mu ,\nu} .\theta ) ) = (f_0(\theta ),\ldots ,f_s(\theta ) )
A_{\beta_{\mu ,\nu}}^{(s)} (\theta)^{-1} ~,
$$
for all $\theta \in {\mathbf T}$. Consequently,
$$
\phi_s ({\mathbf f} \mid_{\lambda} T_{\alpha} )=
$$
$$
\sum_{\mu =1}^r \left(
\begin{array}{c} \vdots \\ ( f_0(\gamma_{\xi_{\mu} (\nu)}
.0),\ldots ,f_s(\gamma_{\xi_{\mu} (\nu)} .0) ) A_{\beta_{\mu
,\nu}}^{(s)} (\gamma_{\xi_{\mu} (\nu)} .0)^{-1}
A_{\alpha_{\mu}}^{(s)} (\gamma_{\nu} .0) \\ \vdots \end{array}
\right)
$$
$$
=  t_{\alpha} \left( \begin{array}{c} \vdots \\
f_0(\gamma_{\nu} .0),\ldots ,f_s(\gamma_{\nu} .0) \\ \vdots
\end{array} \right) = t_{\alpha} \phi_s ({\mathbf f}) ~,
$$
as desired.
\end{proof}

\subsection{} \label{2.5} {\it Proof of Theorem \ref{t1}:} Fix $\lambda \in {\C}$ with $\mathrm{Re}(\lambda) \ge 0$.
Recall (Definition \ref{d2}) that we have fixed a non-negative integer $k$ such that
$$
{\D}_{\lambda ,k}^{\Gamma} = {\D}_{\lambda}^{\Gamma} \cong M_{\lambda}(N) ~,
$$
where the isomorphism is as $\T$-modules, cf. Proposition \ref{p1}. Let $\Phi = (\nu_p \mid p\nmid N)$ be a system of
Hecke eigenvalues occurring in $M_{\lambda}(N)$ and hence also in ${\D}_{\lambda ,k}^{\Gamma}$. Let $0\not =v\in
{\D}_{\lambda ,k}^{\Gamma}$ be a corresponding $\T$-eigenvector, i.e. $T_pv = \nu_p v$ for $p\nmid N$. Recall
(Definition \ref{d2} and Proposition \ref{p2} (1)) that we have a short exact sequence of
$\Delta_1(N)$-modules
$$
0\longrightarrow V_{\lambda ,k} \longrightarrow U_{\lambda ,k} \longrightarrow {\D}_{\lambda ,k} \longrightarrow 0
~,
$$
which gives rise to the very short exact sequence of $\T$-modules
$$
U_{\lambda ,k}^{\Gamma} \stackrel{\alpha}{\longrightarrow} {\D}_{\lambda}^{\Gamma} \stackrel{\beta}{\longrightarrow}
H^1(\Gamma ,V_{\lambda ,k}) ~,
$$
with $\T$ acting as an algebra of cohomological Hecke operators. Suppose that $\beta v=0$ so that $v\in
\mathrm{Im}(\alpha)$. Then $\Phi$ occurs in $\mathrm{Im}(\alpha)$, and because of Theorem \ref{t2} the tautological
homomorphism
$$
\alpha \colon U_{\lambda ,k}^{\Gamma} \longrightarrow \mathrm{Im}(\alpha)
$$
is a surjective homomorphism of $\T$-modules that are finite-dimensional as complex vector spaces. Using Proposition
1.2.2 of \cite{AS}, we may then conclude that $\Phi$ occurs in $U_{\lambda ,k}^{\Gamma}$. According to Theorem \ref{t2}
all eigenvalues $\nu_p$ are then algebraic over the field
$$
F_{\lambda} :={\Q}(\lambda ,\sqrt{n} , n^{\lambda /2} \mid ~n\in \N ) ~.
$$

We then suppose that $\beta v\not = 0$ and will show that $\Phi$ occurs in $H^1(\Gamma ,W_{\lambda ,m})$ for some $0\le
m\le k$. This then finishes the proof of Theorem \ref{t1}.
\smallskip

Recall that according to Proposition \ref{p2} (3) we may consider the spaces $V_{\lambda ,m}$, $0\le m\le k$, as a
filtration
$$
0= V_{\lambda ,0} \le \ldots V_{\lambda ,m-1} \le V_{\lambda ,m} \ldots \le V_{\lambda ,k}
$$
of $\T$-submodules of $V_{\lambda ,k}$ where the successive quotients are isomorphic to the $W_{\lambda ,m}$ as
$\T$-modules. Now, our assumption $\beta v\not = 0$ implies that $\Phi$ occurs in a finite-dimensional sub-$\T$-module
of $H^1(\Gamma ,V_{\lambda ,k})$ namely $\beta ({\D}_{\lambda}^{\Gamma})$; notice, that we must have $k\ge 1$. We shall
assume that $\Phi$ occurs in some finite-dimensional, sub-$\T$-module $X$ of $H^1(\Gamma ,V_{\lambda ,m})$ for some $m$
with $1\le m\le k$ and will show that then either $\Phi$ occurs in $H^1(\Gamma ,W_{\lambda ,m})$ or else in some
finite-dimensional, sub-$\T$-module of $H^1(\Gamma ,V_{\lambda ,m-1})$. By induction on $k$, this gives the desired
conclusion as $V_{\lambda ,0}=0$.
\smallskip

Consider the short exact sequence of $\T$-modules
$$
0\longrightarrow V_{\lambda ,m-1} \longrightarrow V_{\lambda ,m} \longrightarrow W_{\lambda ,m} \longrightarrow
0
$$
coming from Proposition \ref{p2} (3). This gives rise to a long exact sequence of $\T$-modules:
$$
W_{\lambda ,m}^{\Gamma} \longrightarrow H^1(\Gamma ,V_{\lambda ,m-1}) \stackrel{\epsilon}{\longrightarrow} H^1(\Gamma
,V_{\lambda ,m}) \stackrel{\eta}{\longrightarrow} H^1(\Gamma ,W_{\lambda ,m}) ~.
$$
Now, it is easy to see that the space $W_{\lambda ,m}^{\Gamma}$ is 0: Suppose that $f\in W_{\lambda ,m}^{\Gamma}$. As
$f$ is continuous, $f$ is bounded. However, the definition of $j(\gamma ,\theta)$ shows that $j(\gamma_n ,\theta)$,
$n\in \N$, is unbounded if $0\not =\theta \in {\mathbf T}$ where $\gamma_n$ denotes the matrix
$$
\gamma_n := \left( \begin{array}{cc} 1 & n
\\ 0 & 1\end{array} \right) \in \Gamma ~.
$$
So, the definition of the $\Gamma$-structure on $W_{\lambda ,m}$ implies that $f(\theta)$ vanishes for every $\theta
\not =0$ and hence for all $\theta$.

Thus, $\epsilon$ is an injection and the very short exact sequence
$$
\epsilon^{-1} (X) \longrightarrow X\longrightarrow \eta(X)
$$
is an exact sequence of finite-dimensional $\T$-modules. Applying Proposition 1.2.2 of \cite{AS} as above, we conclude
that either $\Phi$ occurs in $\epsilon^{-1} (X)$, and thus a fortiori in $H^1(\Gamma ,V_{\lambda ,m-1})$, or else in
$\eta(X)$ and so also in $H^1(\Gamma ,W_{\lambda ,m})$.
\medskip

\noindent {\bf Remarks:} The reader will notice that our use of continuous coefficients in the above -- as opposed to
$L^2$-coefficients -- is necessitated by the use we made of the evaluation maps $\phi_s$ in the proof of Theorem
\ref{t2} above. Thus, the use of continuous coefficients is indispensable for our approach. However, we wish to remark
here that the there is a certain price to be paid for this, notably the following.

The author does not have concrete examples of triples $(\Gamma, \lambda, s)$ where he can prove that the space
$U^{\Gamma}_{\lambda,s}$ is actually non-zero. However, given the injectivity of the evaluation map $\phi_s$, the proof
of Theorem \ref{t2} shows that any eigenvalue of a Hecke operator $T_{\alpha}$ acting on the space
$U^{\Gamma}_{\lambda,s}$ is also an eigenvalue of the linear operator $t_{\alpha}$ acting on ${\C}^{m(s+1)}$; moreover,
the eigenvalues of $t_{\alpha}$ can -- for any concretely given triple $(\Gamma, \lambda, s)$ -- be computed
numerically. Such numerical experiments seem to indicate that Hecke eigenvalues on the spaces $U^{\Gamma}_{\lambda,s}$
are probably not very interesting, and that one could at the most retrieve packages of Hecke eigenvalues which are
readily recognizable as belonging to certain standard Eisenstein series. Thus, it would appear that the interesting
packages of Hecke eigenvalues should be the ones occurring in the spaces $H^1(\Gamma ,W_{\lambda ,m})$.

There are certain reasons that make it not wholly unreasonable to venture the conjecture that the spaces $H^1(\Gamma
,W_{\lambda ,m})$ are in fact finite-dimensional. For instance, the methods of the papers \cite{BO} and more
specifically \cite{BO1}, might show the way towards analyzing this question. It will be seen however, that because the
$W_{\lambda ,m}$ are spaces of {\it continuous} functions, an attempt to use the methods of these papers to approach
the question of finite-dimensionality of the $H^1(\Gamma ,W_{\lambda ,m})$ will quickly lead to some serious analytical
difficulties.

\end{document}